\documentclass[11pt]{amsart}

\usepackage{fullpage}

\usepackage{url}

\usepackage{amssymb}

\usepackage{cite}

\renewcommand{\a}{\alpha}

\newcommand{\g}{\gamma}
\renewcommand{\d}{\delta}
\newcommand{\D}{\Delta}
\newcommand{\e}{\varepsilon}
\newcommand{\f}{\varphi}
\newcommand{\s}{\sigma}
\renewcommand{\S}{\Sigma}

\renewcommand{\l}{\lambda}

\newcommand{\cF}{{\mathcal F}}
\newcommand{\cC}{{\mathcal C}}

\newcommand{\cT}{{\mathcal T}}

\newcommand{\cD}{{\mathcal D}}
\newcommand{\cR}{{\mathcal R}}

\newcommand{\bR}{\mathbb R}

\newcommand{\bZ}{\mathbb Z}

\newcommand{\bF}{\mathbb F}

\newcommand{\be}{\begin{equation}}
\newcommand{\ee}{\end{equation}}

\renewcommand{\to}{\rightarrow}

\renewcommand{\phi}{\varphi}
\renewcommand{\epsilon}{\varepsilon}

\newcommand{\<}{\langle}
\renewcommand{\>}{\rangle}

\newcommand{\dm}{{\partial M}}

\newcommand{\w}{\widetilde}

\theoremstyle{plain}
\newtheorem{theorem}{Theorem}[section]

\newtheorem{remark}[theorem]{Remark}

\newtheorem{proposition}[theorem]{Proposition}

\theoremstyle{definition}

\def\endproof{\qed \medskip}
\def\blacksquare{\hbox to .60em {\vrule width .60em height .60em}}

\numberwithin{equation}{section}

\begin{document}

\title[ ]{Elliptic regularization of the isometric immersion problem}

\author[ ]{Michael T. Anderson}

\address{Department of Mathematics, Stony Brook University, Stony Brook, N.Y.~11794-3651, USA} 
\email{anderson@math.sunysb.edu}
\urladdr{http://www.math.sunysb.edu/$\sim$anderson}

\thanks{Partially supported by NSF grant DMS 1607479}

\begin{abstract}
We introduce an elliptic regularization of the PDE system representing the isometric immersion of a surface in $\bR^{3}$. The 
regularization is geometric, and has a natural variational interpretation. 
\end{abstract}

\maketitle

\setcounter{section}{0}
\setcounter{equation}{0}

\section{Introduction}

  In this short note, we introduce an elliptic regularization of the equations for isometric immersion of a surface in $\bR^{3}$ 
(or more generally any ambient 3-manifold). Thus we exhibit a smooth curve $\cD_{\e}$, $\e \geq 0$, of differential operators 
which are elliptic for $\e > 0$ for which $\cD_{0}$ is the operator describing isometric immersions. The existence of such a 
regularization is somewhat surprising, since the system of $1^{\rm st}$ order equations for isometric immersions is characteristic 
in all directions and thus seemingly far from elliptic. The regularization $\cD_{\e}$ depends only on geometric data of the immersion. 

\medskip 

  To begin, we recall the global formulation of the problem. Let $\S$ be a closed, orientable 2-dimensional surface, thus a surface of 
genus $g \in \bZ^{\geq 0}$. Any such surface embeds in $\bR^{3}$ and a general immersion 
\be \label{a}
F: \S \to \bR^{3}.
\ee
induces a metric $\g$ on $\S$ by pulling back (or restricting) the Euclidean metric $g_{Eucl}$ to $\S$ via $F$:
\be \label{b}
\g = F^{*}(g_{Eucl}).
\ee
The isometric immersion problem is the converse; given an (abstract) Riemannian metric $\g$ on $\S$, is there 
an immersion $F$ as in \eqref{a} for which \eqref{b} holds? Thus, one is asking which metrics on $\S$ can be ``pictured" as immersed 
surfaces in $\bR^{3}$. A local version of this problem, where $\S$ is replaced by a disc, may be formulated in 
the same way.

   There is a very large literature on this classical problem. This short note is not the place to summarize this in any detail; we refer 
instead to \cite{HH}, \cite{Gr}, \cite{Ho} for background and further references. We only recall that it is a well-known and basic open 
question whether any smooth metric $\g$ on a disc has a neighborhood of $0$ realized by an immersion $F$, i.e.~whether any 
smooth metric locally has a smooth isometric immersion into $\bR^{3}$. Much less is known in general about the global isometric 
immersion problem for compact surfaces. 

  In local coordinates $x^{i}$, $i = 1,2$, on $\S$, the equation \eqref{b} has the form 
\be \label{c}
\sum \partial_{x_{i}}F^{\mu}\cdot \partial_{x_{j}}F^{\mu} = \g_{ij} = \g(\partial_{x_{i}}, \partial_{x_{j}}).
\ee
This is a determined system of three $1^{\rm st}$ order differential equations $1 \leq i \leq j \leq 2$ for three unknown functions $F = 
\{F^{\mu}\}$, $\mu = 1,2,3$. A simple symbol calculation shows that \eqref{c} is not an elliptic system; in fact all directions on the surface are 
characteristic, cf.~\cite{Gr}. As recalled in Section 2, the failure of ellipticity is also easily seen to be a consequence of  Gauss' Theorema Egregium.

   It is well-known that {\em locally}, the isometric immersion system can be reduced to a single scalar equation, the Darboux equation 
$$det D_{\g}^{2}u = K det \g (1 - |\nabla_{\g} u|^{2},$$
for an unknown function $u = u(x,y)$, $|\nabla u|_{\g}^{2} < 1$. Here $K= K_{\g}$ is the Gauss curvature of the immersion $F$, 
and $\nabla_{\g}, D_{\g}^{2}$ are the gradient and Hessian with respect to $\g$. The function $u$ is given locally by 
$u = F\cdot e$, where $e$ is a unit vector in $\bR^{3}$. This Monge-Ampere equation changes type from elliptic when $K > 0$ to 
hyperbolic when $K < 0$ and is degenerate at the locus $K = 0$. The Darboux equation has been the main tool used in understanding 
the local isometric embedding problem, but is not particularly useful for the global isometric embedding or immersion problem.

\medskip 

  To describe the elliptic regularization, given a metric $\g$ on a surface $\S$, let $[\g]$ be the pointwise conformal class of $\g$. 
The full metric $\g$ may then be decomposed into a pair 
\be \label{d}
\g \sim ([\g], \l^{2}),
\ee
where $\l^{2}$ represents the conformal factor with respect to a fixed background metric $\g_{0}$ (for instance of constant curvature) in
the conformal class $\g$. Thus, $\g = \l^{2}\g_{0}$.  Given an immersion $F: \S \to \bR^{3}$, let $\g = F^{*}(g_{Eucl})$ and let $H = H_{F}$ 
denote its mean curvature. 

\begin{theorem} For any $\e > 0$, the data 
\be \label{e}
\cD_{\e}(F) = ([\g], (1-\e)\l^{2} + \e H),
\ee
form a determined elliptic system for an immersion $F: \S \to \bR^{3}$. 
\end{theorem}

  When $\e = 0$, as in \eqref{d} the data $([\g], \l^{2})$ are equivalent to the data $\g$, i.e.~\eqref{e} gives the equations for an 
isometric immersion \eqref{c} when $\e = 0$. Thus one has a smooth path of differential operators $\cD_{\e}$, elliptic for $\e > 0$, with 
$\cD_{0}(F) = \g$ giving the operator for isometric immersions. The local version of Theorem 1.1 is equally valid. 

  In Section 3, we prove that the Fredholm index of $\cD_{\e}$ is zero, for $\S = S^{2}$. This remains unknown for surfaces of 
higher genus, but some partial results on the Fredholm index are discussed in Section 3. We also exhibit a variational (or Lagrangian) 
formulation of the data $\cD_{0}$, $\cD_{1}$ and for data essentially equivalent to $\cD_{\e}$, $\e \in (0,1)$. 

  It would be interesting to approach the isometric immersion problem by studying the behavior of the operators $\cD_{\e}$ as 
$\e \to 0$. For instance, the fact that $\cD_{\e}$ is Fredholm implies that its image is a variety of finite codimension in the target 
space, for all $\e > 0$. It would be of interest to understand the behavior of the kernel (and cokernel) of the linearizations 
$D\cD_{\e}$ as $\e \to 0$, as an approach toward understanding the infinitesimal rigidity of isometric immersions - another 
well-known open problem. We hope to pursue these issues elsewhere.

\section{Elliptic regularization} 
  
  As in the Introduction, let $\S$ be a compact orientable surface. With minor modifications, the discussion below applies equally 
well to the local situation where $\S$ is a disc. 

  Let $Imm^{m+1,\a}(\S)$ be the space of $C^{m+1,\a}$ immersions $F: \S \to \bR^{3}$, where $C^{m+1,\a}$ is the usual 
H\"older space and $m \geq 1$. Let $Met^{m,\a}(\S)$ be the space of $C^{m,\a}$ Riemannian metrics on $\S$. Both of these 
spaces are Banach manifolds and one has a natural map 
\be \label{pid}
\Phi_{D}: Imm^{m+1,\a}(\S) \to Met^{m,\a}(\S),
\ee
$$\Phi_{D}(F) = \g = F^{*}(g_{Eucl}).$$
This is a smooth map of Banach manifolds, since $g_{Eucl}$ is $C^{\infty}$ smooth; cf.~the expression \eqref{c} in 
local coordinates. The isometric embedding problem concerns the description or characterization of the image of $\Phi_{D}$. 

\begin{remark} 
{\rm We recall that Gauss' Theorema Egregium is an obstruction to the map $\Phi_{D}$ being Fredholm, (when $m \geq 2)$. 
(A smooth map $\chi$ is Fredholm if its linearization $D\chi$ at any point is a linear Fredholm map, i.e.~$D\chi$ has finite 
dimensional kernel and cokernel, and is of closed range). 

   Namely, if $\g = \Pi_{D}(F)$, for $F \in Imm^{m+1,\a}(\S)$, then the $2^{\rm nd}$ fundamental form $A$ of $\S$ is in 
$C^{m-1,\a}(\S)$. Let $K$ denote the Gauss curvature of $\g$. Gauss' theorem gives 
$$K = det A,$$
so that $K \in C^{m-1,\a}(\S)$. However the set of metrics $\g \in Met^{m,\a}(\S)$ such that $K \in C^{m-1,\a}(\S)$ is of 
infinite codimension. This contradicts the Fredholm property. 
}
\end{remark}

  Given an immersion $F \in Imm^{m+1,\a}(\S)$, let $[\g]$ be the pointwise conformal class of the induced metric 
$\g$ and let $H$ be the mean curvature of $F$. Let $\cC^{m,\a}(\S)$ be the space of pointwise conformal classes 
of $C^{m,\a}$ metrics on $\S$ and $C^{m-1,\a}(\S)$ the space of $C^{m-1,\a}$ functions $f: \S \to \bR$. 
Define 
\be \label{piH}
\Phi_{H}: Imm^{m+1,\a}(\S) \to \cC^{m,\a}(\S)\times C^{m-1,\a}(\S),
\ee
$$\Phi_{H}(F) = ([\g], H).$$ 
This is a smooth map of Banach manifolds, of mixed Dirichlet-Neumann type (or of mixed intrinsic-extrinsic type). 
Note that $[\g]$ is of first order in $F$, while $H$ is of second order. It is proved in \cite{An2} that the data $([\g], H)$
form an elliptic system for $F$ in the sense of \cite{ADN}. It is worthwhile to include the short proof here.

\begin{proposition}
The data $([\g], H)$ form an elliptic system for $F$. In particular, the map $\Phi_{H}$ in \eqref{piH} is Fredholm.  
\end{proposition}

{\bf Proof:}   The linearization $D\Phi_{H}$ acts on vector fields $X$ along the immersion $F$. Write 
$X = X^{T} + \nu N$, where $X^{T}$ is tangent and $N$ is normal to $\S = Im(F)$. Then 
\be \label{lin1}
\d^{*}X = \d^{*}(X^{T}) + \nu A + d\nu\cdot N,
\ee
so that $(\d^{*}X)^{T} = \d^{*}X|_{T\S} = \d^{*}(X^{T}) + \nu A$. The second term here is lower order in $X$ 
and so does not contribute to the principal symbol. The principal symbol $\s$ of 
the first component of $D\Phi_{H}$ is thus
\be \label{lin2}
\s([(\d^{*}X)^{T}]_{0}) = \s([(\d^{*}(X^{T})]_{0}) = \xi_{i}X_{j} - \frac{\xi_{k}X_{k}}{2}\d_{ij},
\ee
where $X_{i}$, $i = 1,2$ are the components of $X$ tangent to $\S$. Setting this to 0 gives
$$\xi_{1}X_{2} = \xi_{2}X_{1} = 0 \ \ {\rm and} \ \ \xi_{1}X_{1} = \xi_{2}X_{2}.$$
Since $(\xi_{1}, \xi_{2}) \neq (0,0)$, it is elementary to see that the only solution 
of these equations is $X_{1} = X_{2} = 0$. Next, for the mean curvature, one has 
$H_{\d^{*}X}' = -\D \nu - |A|^{2}\nu + X^{T}(H)$, where $\D$ is the Laplacian with respect 
to the induced metric $\g = F^{*}(g_{Eucl})$ and $A$ is the $2^{\rm nd}$ fundamental form of 
$F$.  The leading order symbol acting on $\nu$ is thus $|\xi|^{2}\nu$, which vanishes only if 
$\nu = 0$. Thus, the symbol of $D\Pi_{H}$ is elliptic, so that by the regularity theory for 
elliptic systems, cf.~\cite{ADN}, $D\Pi_{H}$ is Fredholm.

{\endproof}

   In a given conformal class $[\g]$ of metrics on an oriented surface $\S$, consider a metric of constant scalar 
curvature $\g_{0}$. We normalize scalar curvature to $-1$, $0$ or $+1$. In first case, the metric $\g_{0}$ is unique. 
In the second case, $\g_{0}$ is unique up to a scaling, so we assume that $area(\g_{0}) = 1$. In the third case, 
$g_{0}$ is unique up to conformal (M\"obius) transformations of $S^{2}$, i.e.~the action of the conformal group 
$Conf(S^{2})$. In this spherical case, we fix $\g_{0}$ to be the standard round metric, induced by the usual 
embedding of $S^{2}(1) \subset \bR^{3}$. 

  Write then 
$$\g = \l^{2}\g_{0},$$
so that, with the assumptions above, $\l$ is uniquely determined by $\g$. Consider now the data 
\be \label{edata}
\cD_{\e}(F) = ([\g], (1-\e)\l^{2} + \e H),
\ee
for $\e > 0$. The choice $\e = 1$ gives the data \eqref{piH} above while $\e = 0$ gives the ``Dirichlet" data, i.e.~the 
data for an isometric immersion, as in \eqref{a}. 

  Consider then the curve of maps 
\be \label{pie}
\Phi_{\e}: Imm^{m+1,\a}(\S) \to \cC^{m,\a}(\S)\times C^{m-1,\a}(\S),
\ee
$$\Phi_{\e}(F) = ([\g], (1-\e)\l^{2} + \e H).$$ 
This gives a smooth path from conformal/mean curvature data to isometric data. 

\begin{proposition}
For all $\e > 0$, the data \eqref{edata} form an elliptic system for an immersion $F: \S \to \bR^{3}$. The maps 
$\Phi_{\e}$ in \eqref{pie} are smooth Fredholm maps between Banach manifolds.   
\end{proposition}

{\bf Proof:} The proof is exactly the same as the proof of Proposition 2.2. Note that the volume form term $\l^{2}$ is 
of lower differentiability order than $H$. The linearization of the volume form $dv_{\g}$ is determined by 
$tr \d^{*}X = div X = div X^{T} + \nu H$. 

\endproof

 This gives an elliptic regularization of the isometric immersion problem, and so proves Theorem 1.1. 

  The choice of the regularizing term $H$ in \eqref{edata} is of course not unique; it could be replaced for instance by a non-vanishing 
function of $H$; the crucial point in obtaining an elliptic system is to have a scalar function depending on the extrinsic 
geometry of the immersion. 

\begin{remark}
{\rm Propositions 2.2 and 2.3 holds for immersions $F: \S \to (N, g)$ into any complete Riemannian 3-manifold, 
i.e.~the ellipticity of the operator $\Phi_{\e}$, $\e > 0$, is independent of the ambient Riemannian manifold. 

}
\end{remark}

\section{Fredholm index and variational formulation}

  In this section, we compute the Fredholm index of the operators $\cD_{\e}$ for $\e > 0$, at least for $\S = S^{2}$ and 
present initial results for the case of higher genus. We also exhibit a variational interpretation of the operators $\cD_{\e}$, 
(or more precisely, essentially equivalent operators). First, note that the Fredholm index of $\cD_{\e}$, $\e > 0$ is independent of 
$\e$, since the index is deformation invariant. 

\medskip 

  To begin, we recall that the space $Imm^{m+1,\a}(\S)$ is not connected in general; cf.~\cite{Pi} for example for results 
on the number of components of $Imm^{m+1,\a}(\S)$. (The number is of course independent of $(m, \a)$, for $m \geq 1$). 
The Fredholm index of $\cD_{1}$ is constant on each component of $Imm^{m+1,\a}(\S)$ (since the Fredholm 
index is deformation invariant), but the index may apriori vary over different components of $Imm^{m+1,\a}(\S)$. 

  Now recall a famous theorem of Smale \cite{Sm} which states that the space $Imm^{m+1,\a}(S^{2})$ is connected.

\begin{theorem}
The Fredholm index of $\cD_{\e}$, $\e > 0$, on $Imm^{m+1,\a}(S^{2})$ equals $6 = dim(Isom(\bR^{3}))$. 
\end{theorem}

{\bf Proof:} It suffices to compute the index of $\cD_{1}$, with data $([\g], H)$, on the round embedding 
$S^{2}(1) \subset \bR^{3}$. 

  Note that the isometry group $Isom(\bR^{3})$ acts smoothly and freely on the space of immersions $Imm^{m+1,\a}(\S)$ 
via $(F, \iota) \to \iota \circ F$, corresponding to translation, rotation or reflection of $F$. This action 
fixes the target data, i.e.~$\Phi_{1}(\iota \circ F) = \Phi_{1}(F)$. To remove this degeneracy, divide the space 
$Imm$ by this action, and consider the quotient space $Imm_{b}$ of based immersions. 
There is a global slice to this action, i.e.~an inclusion $Imm_{b} \subset Imm$, given by 
fixing a point $p_{0} \in \S$, a unit vector 
$e \in T_{p_{0}}(\S)$ and requiring that $F(p_{0}) = (0,0,0) \in \bR^{3}$,  $T_{p_{0}}(\S) = \bR^{2} 
\subset \bR^{3}$ and with $F_{*}(e) = (a, 0, 0)$ for some $a > 0$. Thus, consider the restricted mapping 
\be \label{Pib}
\w \Phi_{H}: Imm_{b}^{m+1,\a}(\S) \to \cC^{m,\a}(\S) \times C^{m-1,\a}(\S),
\ee
$$\w \Phi_{H}(F) = ([\g], H).$$
Theorem 3.1 then follows from the statement that the Fredholm index of $\w \Phi_{H}$ equals $0$. 

  A variation $X$ of $F$ is non-zero in $T Imm_{b}^{m+1,\a}(\S)$ only if $X$ is not the restriction of a Killing field 
on $\bR^{3}$. As in Proposition 2.2, write $X = X^{T} + \nu N$. Then the variation induced on the data $([\g], H)$ is 
$([\d^{*}X^{T} + \nu A], H'_{\d^{*}X})$. Since the round embedding is umbilic and of constant mean curvature, 
$$([\d^{*}X^{T} + \nu A], H'_{\d^{*}X}) = ([\d^{*}X^{T}], L(\nu)),$$  
where as above $L(\nu) = -\D \nu - |A|^{2}\nu$ is the normal variation of the mean curvature. 

  The kernel thus consists of non-Killing fields such that $\d^{*}(X^{T}) = \f \g$, i.e.~$X^{T}$ is a conformal Killing field 
on $S^{2}(1)$, and functions $\nu$ such that $L(\nu) = 0$. The space of conformal Killing fields on $S^{2}(1)$ is $3$-dimensional. 
Next, since $|A|^{2} = 2$, functions $\nu$ such that $L(\nu) = 0$ are $1^{\rm st}$ eigenfunctions of the 
Laplacian on $S^{2}(1)$. This also forms a $3$-dimensional space, giving a total dimension of $6$. However, a 
$3$-dimensional subspace corresponds to Killing fields (the restriction of infinitesimal translations in $\bR^{3}$ to $S^{2}(1)$). 
It follows that $dim Ker D\w \Phi_{H} = 3$. 

   Regarding the cokernel, variations of data in $Im \w \Phi_{H}$ are of the form $([\d^{*}Y], H'_{\d^{*}Y}) = 
([\d^{*}(Y^{T})], H'_{\nu A})$. Hence the data $(B_{0}, f)$ with $B_{0}$ trace-free is in the cokernel if and only if  
$$\int_{\S}\<\d^{*}(Y^{T}), B_{0}\> + fH'_{\nu A} = 0,$$
for all $Y^{T}$ and $\nu$. Each term must thus vanish separately. Applying the divergence theorem to the first term gives 
$\d B_{0} = 0$ so that $B_{0}$ is transverse-traceless. On $S^{2}$ there are no non-zero such forms (the Teichm\"uller space 
of $S^{2}$ is trivial) so that $B_{0} = 0$. Next $L(f) = 0$ so that $f$ is a $1^{\rm st}$ eigenfunction of the Laplacian, with 
$3$-dimensional eigenspace. Thus $dim Coker D\w \Phi_{H} = 3$ and hence the Fredholm index of $\w \Phi_{H}$ is $0$.

{\endproof}

\begin{remark}
{\rm The index of $\cD_{\e}$ is also independent of the Riemannian metric on $\bR^{3}$. Thus it follows that for any complete 
Riemannian metric $g$ on $\bR^{3}$, $dim Ker \cD_{\e} \geq 6$. In particular, any immersion $F$ realizing a given $([\g], H)$ in 
$(\bR^{3}, g)$ is an element of a $6$-parameter family of immersions realizing such data. 

}
\end{remark}

  It is more difficult to analyse the kernel and cokernel of $D\w \Phi_{H}$ at more general embeddings or for surfaces of higher 
genus. In general, the kernel consists of vector fields $X$ such that 
$$\d^{*}(X^{T}) + \nu A = \f \g\ \ {\rm and} \ \ H'_{\d^{*}X} = 0,$$
for some function $\f: \S \to \bR$.  It is not easy to understand the space of solutions of this system. The cokernel consists of 
pairs $(B_{0}, f)$ such that   
$$\int_{\S}\<Y^{T}, \d B_{0}\> +\nu\<A, B_{0}\> + \nu L(f) + fY^{T}(H) = 0,$$ 
for all $Y^{T}, \nu$. 
Suppose for instance $H = const$. It follows that $B_{0}$ is transverse-traceless and so represents a tangent vector to the 
Teichm\"uller space $T(\cT(\S))$ of $\S$. This gives $dim Coker \geq dim \cT(\S)$. The function $f$ satisfies $L(f) + \<A, B_{0}\> 
= 0$ but again it is difficult to evaluate the dimension of the space of solutions of this equation. 
 
\medskip 

  Next we show that the data \eqref{piH}, \eqref{pid} arise as boundary data for a natural variational problem on the space of metrics 
on a filling of $\S$. This is of independent interest, and gives a new proof and a partial generalization of Theorem 3.1 to higher genus. 

  Let $F: \S \to \bR^{3}$ be an embedding; then $F$ extends to an embedding of a 3-manifold $M$ with $\dm = \S$: 
$\bar F: M \to \bR^{3}$, $\bar F|_{\dm} = F$, and the pull-back $\bar F^{*}(g_{Eucl})$ induces a flat metric on $M$. 
(More generally one may assume that $F$ is an Alexandrov immersion, in that $F$ extends to an immersion 
$\bar F: M \to \bR^{3}$ with $\dm = \S$). This gives a smooth map 
$$\mu: Emb^{m+1,\a}(\S) \to \bF^{m,\a}(M),$$
where $\bF^{m,\a}(M)$ is the Banach manifold of flat metrics on $M$, $C^{m,\a}$ up to $\dm$. Let ${\rm Diff}_{0}^{m+1,\a}(M)$ be the 
group of $C^{m+1,\a}$ diffeomorphisms of $M$ which equal the identity on $\dm = \S$. This acts freely and smoothly on 
$\bF^{m,\a}(M)$ and let $\cF^{m,\a}(M)$ be the quotient space (the moduli space of flat metrics on $M$). The map $\mu$ induces a 
smooth map 
$$\mu: Emb^{m+1,\a}(\S) \to \cF^{m,\a}(M).$$

  Now consider the smooth map 
$$\Pi_{H}: \cF^{m,\a}(M) \to \cC^{m,\a}(\S)\times C^{m-1,\a}(\S),$$
\be \label{PIH}
\Pi_{H}(g) = ([\g], H).
\ee
Note that $\Pi_{H}\circ \mu = \Phi_{H}$ (when $\Pi_{H}$ is restricted to flat embeddings). Similarly, one has a smooth map 
$$\Pi_{D}: \cF^{m,\a}(M) \to Met^{m,\a}(\S),$$
\be \label{PID}
\Pi_{D}(g) = \g,
\ee
with $\Pi_{D}\circ \mu = \Phi_{D}$. 

   To describe the variational formulation, let $Met^{m,\a}(M)$ be the Banach space of $C^{m,\a}$ metrics on $M$. 
Consider first the well-known Einstein-Hilbert action with Gibbons-Hawking-York boundary term \cite{GH}, \cite{Y}.  
Thus let 
$$I_{D}: Met^{m,\a}(M) \to \bR,$$
\be \label{var}
I_{D}(g) = \int_{M}R_{g}dV_{g} + 2\int_{\S}H dv_{\g},
\ee
where $R_{g}$ is the scalar curvature of $g$ and $H$ is the mean curvature of $g$ at $\dm = \S$ with respect to the outward 
normal $N$. The linearization of the scalar curvature $R_{g}$ in the direction $h$ is given by $R'(h) = -\D tr h + \d \d h - \<Ric, h\>$, 
while in geodesic normal coordinates near $\dm = \S$, $H'_{h} = \frac{1}{2}N(tr h)$. From this, a straightforward computation using 
the divergence theorem shows that the linearization of $I$ at $g$ is given by 
\be \label{dvar}
d(I_{D})_{g}(h) = -\int_{M}\<E_{g}, h\>dV_{g} -\int_{\dm}\<\tau, h^{T}\>dv_{\g},
\ee
where $E_{g} = Ric_{g} - \frac{R_{g}}{2}g$ is the Einstein tensor of $(M, g)$, $h^{T} = h|_{\dm}$ and 
$\tau = A - H\g$ is the conjugate momentum (with respect to the functional $I_{D}$); the expression \eqref{dvar} is 
well-known, cf.~\cite{GH}, \cite{Y}. In particular, \eqref{dvar} shows that critical points of $I_{D}$ on the space of metrics 
with fixed boundary metric $\g$ on $\dm$ (Dirichlet data on $\dm$) are flat metrics.

  Essentially the same computation shows that the data $([\g], H)$ also arise as boundary data of a natural variational problem, 
using a slight modification of the Gibbons-Hawking-York boundary term. Thus as in \eqref{var}, let 
$$I_{H}: Met^{m,\a}(M) \to \bR,$$
\be \label{var1}
I_{H}(g) = \int_{M}R_{g}dV_{g} + \int_{\S}H dv_{\g}.
\ee
A similar computation as above (cf.~\cite{An2}) gives 
\be \label{dvar1}
d(I_{H})_{g}(h) = -\int_{M}\<E_{g}, h\>dV_{g} -\int_{\S}(\<A_{0}, h_{0}\> + H'_{h})dv_{\g},
\ee
where $A_{0}, h_{0}$ are the trace-free parts of $A$ and $h^{T}$ respectively (with respect to $\g$). 
In particular, writing $H'_{h}$ as $H (\log H)'_{h}$ in \eqref{dvar1}, one sees that $(A_{0}, (\log H)')$ are dual (conjugate) to the data 
$([\g], H)$ with respect to $I_{H}$. As before, critical points of $I_{H}$ on the space of metrics with fixed conformal class and 
mean curvature are flat metrics. 
 
  The $2^{\rm nd}$ variation of either of the functionals $I_{D}$ or $I_{H}$ leads to ``self-adjoint" properties of the boundary and 
bulk terms. Thus, let $h$, $k$ be a pair of infinitesimal flat deformations of a flat metric $g$ on $M$ and let $g_{s,t} = g + th + sk$. 
Using the equality of mixed $2^{\rm nd}$ derivatives 
$$\frac{\partial^{2}I_{D}}{\partial s \partial t} = \frac{\partial^{2}I_{D}}{\partial t \partial s}$$ 
leads to the relations   
$$\int_{\dm}\<\tau'_{k} + a_{D}(k), h^{T}\>dv_{\g} = \int_{\dm}\<\tau'_{h} +a_{D}(h), k^{T}\>dv_{\g},$$
where $a_{D}(k) = -2\tau \circ k + \frac{1}{2}(tr_{\g}k)\tau$ arises from the variation of the metric and volume form. Similarly 
the $2^{\rm nd}$ variation of $I_{H}$ on $h$, $k$ gives 
$$\int_{\S}(\<(A_{0})'_{k} + a_{H}(k), h_{0}\>  + \tfrac{1}{2} H'_{h} tr k )dv_{\g} = 
\int_{\S}(\<(A_{0})'_{h} + a_{H}(h), k_{0}\> + \tfrac{1}{2} H'_{k} tr h) dv_{\g},$$
where $a_{H}(k) = -2A_{0}\circ k + \frac{1}{2}(tr_{\g}k)A_{0}$. (Here we have also used the fact that 
$\partial^{2}H/\partial s \partial t = \partial^{2}H/\partial t \partial s$). 

  Similarly, if $h$ and $k$ are any variations of the metric $g$, so $h, k \in T_{g}Met^{m,\a}(M)$ with either 
\be \label{dbc}
h^{T} = k^{T} = 0,
\ee
at $\dm$ (for the functional $I_{D}$) or 
\be \label{hbc}
([h^{T}], H'_{h}) = ([k^{T}], H'_{k} ) = 0,
\ee
at $\dm$ (for the functional $I_{H}$) then 
\be \label{bsa}
\int_{M}\<E'_{h}, h\> = \int_{M}\<h, E'_{k}\>.
\ee
This has the following essentially standard consequence. 

\begin{theorem}
The map $\Pi_{H}$ is Fredholm, of Fredholm index 0. 
\end{theorem}

{\bf Proof:}  Given a (background) flat metric $\w g \in \bF^{m,\a}$, we work in the divergence-free gauge with respect to $\w g$ 
for the action of  ${\rm Diff}_{0}^{m+1,\a}(M)$ on $Met^{m,\a}(M)$. Thus consider the divergence-gauged Einstein operator 
$$\Phi(g) = E(g) + \d^{*}\d_{\w g}(g),$$ 
and its linearization at $g = \w g$, 
$$2L(h) = E'(h) + \d^{*}\d h.$$
This is an elliptic operator and boundary conditions $(\d h, [h], H'_{h})$ form an elliptic boundary value problem for $L$; see 
\cite{An1} for a proof. Moreover, if $\d h = 0$ on $\dm$, then $L(h) = 0$ implies $E'(h) = 0$ and $\d h = 0$. 

   Let $S_{0}^{m,\a}(M) \subset S^{m,\a}(M) = TMet^{m,\a}(M)$ be the subspace of variations $h$ of $g$ such that 
\be \label{bc1}
\d h = 0, \ \ [h^{T}] = 0, \ \ H'_{h} = 0,
\ee
at $\dm$. It follows from \eqref{hbc} and \eqref{bsa} that 
$$L: S_{0}^{m,\a}(M) \to S^{m-2,\a}(M),$$
is a formally self-adjoint elliptic operator. By the Fredholm alternative, 
\be \label{fred}
Im L \oplus K = S^{m-2,\a}(M),
\ee
where $K = Ker L \subset S_{0}^{m+1,\a}(M)$. This gives a natural identification of the kernel and cokernel of $L$ and in 
particular, 
$$dim Ker L = dim Coker L.$$
Now given (arbitrary) boundary data $([h^{T}], H') \in T(\cC^{m,\a}(\S)\times C^{m-1,\a}(\S))$ with $\d h = 0$ at $\S$, let $h$ be a 
$C^{m,\a}(M)$ extension of the boundary data $([h^{T}], H')$ to $M$. Let $z = L(h)$ and via \eqref{fred} decompose $z$ uniquely as 
$z = L(\bar h) + k$ with $k \in Ker L$. Then $\w h = h - \bar h$ satisfies $L(\w h) = k$ and the boundary data of $\w h$ are given 
by $([h^{T}], H')$. 

  This shows that there is subspace of codimension equal to $dim K$ in the space of boundary data for which there is an extension
$\w h$ such that $L(\w h) = 0$ and hence $E'(\w h) = 0$. It follows that the mapping \eqref{PIH} is Fredholm and of Fredholm index 0. 

{\endproof}

\begin{remark}
{\rm By Smale's theorem \cite{Sm}, Theorem 3.3 implies Theorem 3.1. However, Theorem 3.3 does not give an immediate generalization 
of Theorem 3.1 to surfaces of genus $g \geq 1$ since the flat metrics in $\cF^{m,\a}$ may have non-trivial holonomy. The flat 
deformations $h$ in \eqref{PIH} or \eqref{PID} are {\it locally} of the form $h = \d^{*}X$ for a vector field $X$, but not necessarily 
globally of this form. The quotient space of $\cF^{m,\a}$ modulo the action of the full group ${\rm Diff}^{m+1,\a}(M)$ of $C^{m+1,\a}$ 
diffeomorphisms $\f: M \to M$ is the representation variety $\cR(M)$; the space of group homomorphisms $\pi_{1}(M) \to Isom(\bR^{3})$, 
cf.~\cite{L}, \cite{T}, \cite{HK}. 

  Thus, flat deformations of immersions  $\S \to \bR^{3}$ typically include a deformation of the holonomy (from trivial to non-trivial). 
An explicit example of this behavior for $\S = T^{2}$ is exhibited in \cite{HL}. 
}
\end{remark}

  Finally, we show that a modification for the data $\cD_{\e}$ in \eqref{edata} also have a variational interpretation for 
$\e \in [0,1]$. Consider then the linear combination 
$$I_{\e}(g) = (1-\e)I_{D} + \e I_{H} = \int_{M}R dV_{g} + (2-\e)\int_{\dm}H dv_{\g}.$$
One has $dv_{\g} = \l^{2}dv_{\g_{0}}$, so that $tr h / 2 = (dv_{\g})'/dv_{\g} = 2\l'/\l$. The variational derivative of $I_{\e}(h)$ is 
thus given by 
$$d(I_{\e})_{g}(h) = -\int_{M}\<E, h\> dV_{g} - \int_{\S}(\<(1-\e)\tau + \e A_{0}, h\> + \e H'_{h})dv_{\g}.$$
Since $\tau = A - H\g$, one computes $\<(1-\e)\tau + \e A_{0}, h\> + \e H'_{h} = \<A_{0}, h\>  - (1-\e)\tfrac{H}{2}\<\g, \tfrac{tr h}{2}\g\> 
+ \e H'_{h} = \<A_{0}, h\> - (1-\e)\tfrac{H}{2} tr h + \e H'_{h}$, so that 
$$d(I_{\e})_{g}(h)  = -\int_{M}\<E, h\> dV_{g} - \int_{\S}(\<A_{0}, h\> - (1-\e)\tfrac{H}{2} tr h + \e H'_{h})dv_{\g}.$$
Also, $(1-\e)\tfrac{H}{2} tr h - \e H' = H((1-\e)2\l'/\l - \e H'_{h}/H) = H((\log \l^{2(1-\e)}H^{-\e})'$. Thus critical points of $I_{\e}$ 
on metrics with fixed conformal class (so $h_{0} = 0$) and with fixed product $\l^{2(1-\e)}H^{-\e}$ are exactly the flat metrics. 
In particular, the data 
$$\w \cD_{\e}(F) = ([\g], \l^{2(1-\e)}H^{-\e}),$$
have a natural variational formulation analogous to that of $\cD_{D}$ or $\cD_{H}$. The same proof as that in 
Proposition 2.2 shows that this data is elliptic. Moreover, replacing the definition of $S_{0}^{m,\a}(M)$ in \eqref{bc1} by the 
boundary conditions 
$$(\d h, [h], (\l^{2(1-\e)}H^{-\e})'_{h}) = (0,0,0),$$
the same proof as that in Theorem 3.3 shows that the map 
$$\Pi_{\e}: \cF^{m,\a}(M) \to \cC^{m,\a}(\S)\times C^{m-1,\a}(\S),$$
\be \label{PIE}
\Pi_{\e}(g) = ([\g], \l^{2(1-\e)}H^{-\e}),
\ee
is Fredholm, of Fredholm index 0.

 With the boundary conditions defining $S_{0}^{m,\a}(M)$ in \eqref{bc1} replaced by the Dirichlet boundary conditions 
$\d h = 0$ and $h^{T} = 0$ at $\dm$, the operator $L$ and so $E'$ is still formally self-adjoint. However, it is no longer elliptic. 
It would be interesting to understand the behavior of the kernel and cokernel of $L$ on these spaces as $\e \to 0$.

\bibliographystyle{plain}

\end{document}